\documentclass[11pt]{article}
\usepackage[T1]{fontenc}
\usepackage{lmodern}
\usepackage[margin=1in]{geometry}
\usepackage{amsmath,amssymb,amsthm,mathtools}
\usepackage{microtype,booktabs,enumitem}
\usepackage{xcolor}
\usepackage[colorlinks=true,linkcolor=blue!45!black,citecolor=blue!45!black,urlcolor=blue!45!black]{hyperref}
\hypersetup{pdftitle={Isospectral Density Design and Distinct Global Shape Transitions in Weighted Yukawa Droplets},pdfauthor={Qiushuang Liu and Ning Zhang},pdfsubject={A constrained inverse design theorem for a nonlocal isoperimetric problem}}
\numberwithin{equation}{section}
\newtheorem{theorem}{Theorem}[section]
\newtheorem{lemma}[theorem]{Lemma}
\newtheorem{proposition}[theorem]{Proposition}

\theoremstyle{remark}
\newcommand{\R}{\mathbb R}
\newcommand{\Sph}{\mathbb S^2}
\newcommand{\dd}{\,\mathrm d}
\newcommand{\one}{\mathbf 1}
\newcommand{\N}{\mathcal N}
\newcommand{\J}{\mathcal J}
\newcommand{\ip}[2]{\langle #1,#2\rangle}
\newcommand{\norm}[1]{\lVert #1\rVert}
\newcommand{\SO}{\mathrm{SO}(3)}
\DeclareMathOperator{\tr}{tr}

\setlist[enumerate]{leftmargin=*,itemsep=3pt,topsep=5pt}
\title{Isospectral Density Design and Distinct Global\\[3pt] Shape Transitions in Weighted Yukawa Droplets}
\author{%
Qiushuang Liu \qquad Ning Zhang\\[5pt]
{\small School of Mathematics and Statistics}\\
{\small Huazhong University of Science and Technology}\\
{\small 1037 Luoyu Road, Wuhan, Hubei 430074, China}\\[3pt]
{\small \href{mailto:d202580018@hust.edu.cn}{\texttt{d202580018@hust.edu.cn}}\qquad
\href{mailto:nzhang2@hust.edu.cn}{\texttt{nzhang2@hust.edu.cn}}}}
\date{}
\begin{document}
\maketitle
\begin{abstract}
We study Euclidean perimeter plus Yukawa repulsion with a fixed radial
source density and a fixed ordinary volume. For each sufficiently small
positive screening parameter, we construct two smooth, positive, strictly
increasing densities whose constrained shape Hessians at the unit ball
agree for every coupling. Their source masses, Yukawa moments, and ball
self-energies also agree. Just above their common first global transition
threshold, all minimizers form one rotation orbit of tetrahedral shapes
for one density and one orbit of sectoral threefold shapes for the other.
Both transitions are continuous, and the minimizing shapes are strictly
convex. We obtain the two densities by fixing the boundary data that
determine the Hessian and changing the higher derivatives that determine
the reduced quartic energy. We compute a positive quartic coefficient,
restore the three scalar observables by interior corrections, and add a
finite exterior density barrier to control all finite-perimeter competitors.
\end{abstract}

\section{The question and the main result}\label{sec:intro}

The shape Hessian at a sphere determines its unstable modes. It need not
determine the shapes that minimize the nonlinear energy. We construct two
radial source densities with the same ball Hessian, the same source mass,
Yukawa moment, and ball self-energy, and the same first global transition
threshold. The two models have different minimizing shapes just above
that threshold.

For $\kappa>0$, let
\begin{equation}\label{eq:model}
 G_\kappa(x)=\frac{e^{-\kappa|x|}}{4\pi|x|},\qquad
 F_{\beta,\rho}(E)=P(E)+\beta I_\rho(E),\qquad
 I_\rho(E)=\frac12\iint_{E\times E}\rho(|x|)G_\kappa(x-y)\rho(|y|)\dd x\dd y.
\end{equation}
Here $P$ is Euclidean perimeter, $\beta\ge0$, and $\rho$ is fixed in space.
The admissible sets have finite perimeter and ordinary volume $|E|=|B|$,
where $B=B_1(0)$. Sets are identified up to null sets. The integral is allowed
to be $+\infty$; such sets are irrelevant for minimization when $\beta>0$.
At $\beta=0$ the energy means perimeter alone.
The density weights both sources in the interaction. The perimeter and
volume are unweighted.

Write $V_\ell$ for the real spherical harmonics of degree $\ell$, with the
$L^2(\Sph,\dd\omega)$ inner product. The constrained Hessian at $B$ is
diagonal on $\bigoplus_{\ell\ge1}V_\ell$. We use ``isospectral'' to mean
equality of the Hessian operators on each degree. Define
\begin{equation}\label{eq:observables}
 M_0(\rho)=4\pi\int_0^1r^2\rho(r)\dd r,\qquad
 M_\kappa(\rho)=\int_0^1r\sinh(\kappa r)\rho(r)\dd r,
\end{equation}
and let $\beta_g(\rho)$ be the supremum of the positive couplings at which
$B$ is a global minimizer. The construction below proves that this set of
couplings is exactly $(0,\beta_g]$.

The two unit critical harmonics used throughout are
\begin{equation}\label{eq:representatives}
 e_T=\sqrt{\frac{105}{4\pi}}xyz,\qquad
 e_S=\sqrt{\frac{35}{32\pi}}(x^3-3xy^2)\quad\text{on }\Sph.
\end{equation}
Their full orthogonal stabilizers are the tetrahedral group $T_d$ and the
sectoral threefold group $D_{3h}$, respectively. The associated droplets
are three-dimensional. The term ``sectoral'' describes the spherical harmonic.

\begin{theorem}[Identical linear and scalar data, distinct first global shapes]
\label{thm:main}
There exists $\kappa_0>0$ such that, for every fixed
$0<\kappa<\kappa_0$, there are two radial densities $\rho_T,\rho_S$,
smooth as functions on $\R^3$, and numbers $\beta_*,\varepsilon>0$ with
the following properties.
\begin{enumerate}
\item The densities are positive, satisfy $\rho_i'(r)>0$ for $r>0$,
coincide outside a compact annulus, and have quadratic growth at infinity.
Their constrained ball Hessians agree for every $\beta\ge0$, and
\[
 M_0(\rho_T)=M_0(\rho_S),\quad
 M_\kappa(\rho_T)=M_\kappa(\rho_S),\quad
 I_{\rho_T}(B)=I_{\rho_S}(B).
\]
\item For each density a global minimizer exists for every $\beta>0$.
The centered ball is the unique global minimizer for
$0<\beta\le\beta_*$. Its Hessian is positive definite for
$0<\beta<\beta_*$ and has kernel exactly $V_3$ at $\beta_*$.
Also,
\[
 \beta_g(\rho_T)=\beta_g(\rho_S)=\beta_*,\qquad
 \lim_{\kappa\downarrow0}\beta_*(\kappa)=175.
\]
\item For $\beta_*<\beta<\beta_*+\varepsilon$, all global minimizers
of the $i$th model form exactly one $\SO$ orbit of nonspherical,
strictly convex $C^{2,\gamma}$ domains, for any fixed $0<\gamma<1$.
Their full orthogonal stabilizers are $T_d$ for $i=T$ and $D_{3h}$ for
$i=S$.
\item Put $\mu=\beta-\beta_*$. Up to rotation, the boundary is
$r_i(\omega)^3=1+3w_i(\omega)$ with
\begin{equation}\label{eq:main-asymptotic}
 w_i=t_i e_i+O_{C^{2,\gamma}}(t_i^2),\qquad
 t_i^2=\frac{\alpha\mu}{4a_i}+O(\mu^2),\qquad
 F_{\beta,\rho_i}(E_i)-F_{\beta,\rho_i}(B)
 =-\frac{\alpha^2\mu^2}{16a_i}+O(\mu^3),
\end{equation}
where $\alpha,a_T,a_S$ are positive. The first global transition is
continuous, with different leading angular shapes.
\end{enumerate}
\end{theorem}

The densities depend on the fixed screening parameter and the unit
reference volume. We do not require one pair to have the same Hessians at
all screening parameters or radii. The theorem compares two functionals;
it makes no claim that both shape types minimize the same functional.
The boundary-layer width, exterior barrier strength, and transition
interval may depend on $\kappa$; no uniform bounds on these choices as
$\kappa\downarrow0$ are claimed.

\paragraph{Relation to earlier work.}
Frank \cite{Frank} constructed nonspherical volume-constrained equilibria
of the Coulomb liquid drop model through bifurcation from a ball.
De Regibus, Grossi, and Musso \cite{DRGM} study further bifurcations for
nonlocal isoperimetric problems. For screened interactions, Bronsard,
Merlet, and Pegon \cite{BMP} establish connected nonball minimizers in
Yukawa and truncated Coulomb models, with nonradial minimizers for the
latter. Bronsard, DeMason, and Topaloglu \cite{BDT} study existence and
sharp ball stability for the unweighted Yukawa problem. These papers concern bifurcation or minimization for prescribed models.
Here we choose the source density while fixing the ball Hessian and three
scalar observables.

Weights and confinement have also been studied in related liquid drop
problems. Alama, Bronsard, Topaloglu, and Zuniga \cite{ABTZ} place the
density in the perimeter, while Alama, Bronsard, Choksi, and Topaloglu
\cite{ABCT} use a background potential. Their models differ from
\eqref{eq:model}, whose ordinary perimeter and volume are fixed while a
radial source density multiplies both factors of the interaction.

The $O(3)$ bifurcation theory, the two orbit types, and the quartic
invariants of harmonic cubics are classical \cite{IG,CLM,CHQZ,BvV}.
We realize the two orbit types using positive increasing Yukawa source
densities with the same spectral and scalar data. We also prove that the
first global thresholds agree and classify all minimizers just above
them. We include the invariant calculations with the normalizations
used in the proof.

\paragraph{Outline of the proof.}
The ball Hessian depends only on the boundary value, the first boundary
derivative, and one interior Yukawa moment. We change higher boundary
derivatives in a thin layer and use three interior corrections to restore
the observables. We choose the second derivative so that the radial
quartic coefficient is positive, then choose the third derivative to
select either angular orbit. A common finite increase of the exterior
density excludes distant competitors and leaves all near-ball energies
unchanged.

\section{Linear data and exact density realization}\label{sec:density}

\subsection{The complete ball Hessian}
Let $u=G_\kappa*(\rho\one_B)$, $q=\rho(1)$, $p=\rho'(1)$, and
$U=u(1)$. The radial Green formula gives
\begin{equation}\label{eq:radial}
 u(r)=\frac{e^{-\kappa r}}{\kappa r}M_\kappa(\rho)\quad(r\ge1),
 \qquad U=\frac{e^{-\kappa}}\kappa M_\kappa(\rho),\qquad
 u'(1)=-(\kappa+1)U.
\end{equation}
Here $u'(1)$ differentiates the fixed ball potential in space.

\begin{lemma}[Spectrum and its finite data]\label{lem:spectrum}
The ball is volume-constrained critical for every $\beta$, and the
Hessian eigenvalue on $V_\ell$, with multiplicity $2\ell+1$, is
\begin{equation}\label{eq:spectrum}
 \lambda_\ell(\beta)=\ell(\ell+1)-2+
 \beta\{q^2s_\ell+C_\rho\},\qquad
 C_\rho=U[p-(\kappa+1)q],\qquad \ell\ge1,
\end{equation}
where $s_\ell=I_{\ell+1/2}(\kappa)K_{\ell+1/2}(\kappa)$.
For fixed $\kappa$, $s_\ell$ decreases strictly to zero, and
\begin{equation}\label{eq:sl-bound}
 \frac1{2\ell+1+\kappa^2[(2\ell+3)^{-1}+(2\ell-1)^{-1}]}
 \le s_\ell\le\frac1{2\ell+1}\quad(\ell\ge1).
\end{equation}
The data $q,p,M_\kappa$ determine the whole Hessian for every coupling.
\end{lemma}
\begin{proof}
The first variations are $DP[h]=\int Hh$, $DI[h]=\int\rho u_Eh$,
and $D|E|[h]=\int h$. At the sphere the multiplier is $2+\beta qU$.
The second variation of the Lagrangian is
\[
 \int_{\Sph}(|\nabla h|^2-2h^2)
 +\beta C_\rho\int_{\Sph}h^2
 +\beta q^2\iint_{\Sph\times\Sph}G_\kappa(\omega-\sigma)
 h(\omega)h(\sigma)\dd\omega\dd\sigma.
\]
The local term is $\partial_r(\rho u)(1)$; the remaining term is the
single-layer operator. On $V_\ell$ its eigenvalue is the reciprocal of
the jump between the interior and exterior logarithmic radial
derivatives. The normalized solutions are proportional to
$r^{-1/2}I_{\ell+1/2}(\kappa r)$ and
$r^{-1/2}K_{\ell+1/2}(\kappa r)$, respectively. Their Wronskian gives
$s_\ell$ as stated.

The same eigenvalue has the variational formula: $s_\ell^{-1}$ is the minimum over $f(1)=1$ of
\[
 \int_0^\infty\{r^2|f'|^2+\ell(\ell+1)|f|^2+\kappa^2r^2|f|^2\}\dd r.
\]
Comparing degrees at the higher-degree minimizer proves strict
monotonicity. Completing squares on $(0,1)$ and $(1,\infty)$ gives
the lower energy bound $2\ell+1$. Testing with $r^\ell$ inside and
$r^{-\ell-1}$ outside gives the upper bound on the minimum and proves
\eqref{eq:sl-bound}. These trial functions have finite massive energy
for $\ell\ge1$.
\end{proof}

At a known nonzero coupling, two distinct labeled
eigenvalues determine $q^2$ and $C_\rho$: if
$T_\ell=(\lambda_\ell-\ell(\ell+1)+2)/\beta$, then
$q^2=(T_\ell-T_j)/(s_\ell-s_j)$ and $C_\rho=T_\ell-q^2s_\ell$.
We preserve the stronger data $q,p,M_\kappa$.
The translation modes $V_1$ are retained, since the density is fixed in
space.

\begin{lemma}[A common degree-three crossing]\label{lem:crossing}
Set
\begin{equation}\label{eq:baseline}
 \rho_0(r)=a_0+b_0r^2,\qquad
 a_0=\frac{\sqrt7-1}{2},\qquad b_0=\frac{3-\sqrt7}{2}.
\end{equation}
For all sufficiently small positive $\kappa$, this density has a first
positive Hessian crossing at $\beta_*(\kappa)$ with kernel exactly
$V_3$. At the crossing the complement is uniformly coercive in $H^1$,
and $\beta_*(\kappa)\to175$.
\end{lemma}
\begin{proof}
At $\kappa=0$, $q=1$, $p=3-\sqrt7$, $U=(2+\sqrt7)/15$, and
$C_\rho=U(p-1)=-1/5$. This gives
\[
 \lambda_\ell(\beta)=\ell(\ell+1)-2+
 \beta\left(\frac1{2\ell+1}-\frac15\right),\qquad
 \lambda_\ell(175)=\frac{(\ell-3)(2\ell^2+9\ell-46)}{2\ell+1}.
\]
For integer $\ell\ge1$ this is zero only at $\ell=3$ and positive
otherwise. The degree-three slope is $-2/35$, while the degree-one
slope is $2/15>0$. At small $\kappa>0$ set
\begin{equation}\label{eq:alpha}
 \alpha=-(s_3+C_{\rho_0})>0,\qquad \beta_*=10/\alpha.
\end{equation}
Continuity preserves the positive gap for every finite set of
noncritical degrees. For the remaining degrees, $s_\ell\ge0$ and
$C_{\rho_0}$ and $\beta_*$ are bounded, so
$\lambda_\ell(\beta_*)\ge\ell(\ell+1)-C$ uniformly. This proves the $H^1$ gap for all noncritical degrees. Since every eigenvalue is affine in $\beta$, positivity at
$\beta=0$ and $\beta=\beta_*$, together with the positive degree-one
slope, proves positivity for $0<\beta<\beta_*$. The degree-three
eigenvalue is exactly $-\alpha(\beta-\beta_*)$.
\end{proof}

\subsection{Changing finite boundary derivatives with exact compensation}

\begin{lemma}[Finite-jet realization]\label{lem:jets}
Fix $\kappa>0$. Given arbitrary finite values for
$\rho''(1),\rho'''(1)$, there is a smooth positive radial density
$\rho$ realizing them and preserving
\[
 \rho(1),\quad \rho'(1),\quad M_0,\quad M_\kappa,\quad I(B)
\]
of the baseline \eqref{eq:baseline}. Its difference from the baseline
has compact support away from the origin, and can be arbitrarily small
in $C^1$. The density satisfies $\rho'(r)>0$ for every $r>0$.
\end{lemma}
\begin{proof}
With a cutoff $\chi=1$ near zero, first add
\[
 d_\epsilon(r)=\chi((r-1)/\epsilon)
 \left\{\frac{a}{2}(r-1)^2+\frac{b}{6}(r-1)^3\right\}.
\]
The constants $a,b$ are the prescribed increments. This changes the second
and third derivatives, leaves the value and first derivative fixed, and tends to zero
in $C^1$ as $\epsilon\to0$ for fixed $a,b$.

Choose an interval $J\Subset(0,1)$ disjoint from this layer. On
$C_c^\infty(J)$ consider the derivatives of the three observables:
\[
 L_0(d)=4\pi\int r^2d,\qquad
 L_Y(d)=\int r\sinh(\kappa r)d,\qquad
 L_u(d)=4\pi\int r^2u_0d,
\]
where $u_0=G_\kappa*(\rho_0\one_B)$. These functionals are independent.
The first two weights are independent because
$\sinh(\kappa r)/r$ is strictly increasing. Dependence of the third
would imply $u_0=A+B\sinh(\kappa r)/r$ on $J$; applying
$-\Delta+\kappa^2$ would make $\rho_0$ constant there, a contradiction.
Choose three test functions $d_1,d_2,d_3$ with invertible derivative
matrix. The finite-dimensional implicit function theorem applied to
\[
 (c_1,c_2,c_3)\longmapsto
 \big(M_0,M_\kappa,I(B)\big)
 \left(\rho_0+d_\epsilon+\sum c_jd_j\right)
\]
gives small corrections restoring all three values exactly.
We apply the theorem with the density in $C^0$ on the fixed ball. The first
two functionals are linear, and $I(B)$ is a continuous quadratic form.
The higher derivatives of $d_\epsilon$ need not be small. The corrections
are $C^1$ small and supported away from the boundary, so they leave
the prescribed derivatives unchanged. Positivity and strict increase
follow from the positive minima of $\rho_0$ and $\rho_0'$ on the compact
perturbation support. Outside that support the baseline is unchanged.
\end{proof}

The total perturbation preserves $M_\kappa$. An individual interior
correction may have nonzero moment to compensate for the boundary layer.
If an interior perturbation $d$ supported in $B_a$ has zero moment, then
$G_\kappa*d=0$ outside $B_a$ by \eqref{eq:radial}. For every $E\supset B_a$,
the resulting change in $I(E)$ is independent of $E$, so relative shape
energies are unchanged.

\section{Exact reduction and quartic control}\label{sec:reduction}

\subsection{Reduction in a strong space, comparison in a weaker norm}
Use the exact volume coordinate
\begin{equation}\label{eq:volume}
 r_w^3=1+3w,\qquad \int_{\Sph}w=0,\qquad
 E_w=\{r\omega:0\le r<r_w(\omega)\}.
\end{equation}
Then $|E_w|=|B|$, and $r_w-1=w+O(w^2)$. Fix $0<\gamma<1$,
put $V=V_3$, and let $W$ be its $L^2$ orthogonal complement among
mean-zero functions. For $\mu=\beta-\beta_*$ define
$\mathcal E(w,\mu)=F_{\beta_*+\mu,\rho}(E_w)-F_{\beta_*+\mu,\rho}(B)$.

\begin{lemma}[Exact reduction and energy comparison]\label{lem:reduction}
For every fixed smooth radial density with the ball data of
Lemma~\ref{lem:crossing}, there is a smooth, $O(3)$-equivariant map
$z:V\times\R\to W\cap C^{2,\gamma}$ near $(0,0)$ such that
\[
 z(0,\mu)=0,\quad D_vz(0,\mu)=0,\quad
 P_W\nabla\mathcal E(v+z(v,\mu),\mu)=0.
\]
The reduced energy $\Psi(v,\mu)=\mathcal E(v+z(v,\mu),\mu)$ is even
in $v$ and has the expansion
\begin{equation}\label{eq:reduced}
 \Psi(v,\mu)=-\frac\alpha2\mu\norm v_2^2+Q_\rho(v)
 +O(\norm v_2^6+|\mu|\norm v_2^4),
\end{equation}
with the corresponding bounds after up to two $v$ derivatives.
For every sufficiently small mean-zero $C^1$ graph, write
$w=v+z(v,\mu)+y$, $y\in W$. There are $c,C>0$ such that
\begin{equation}\label{eq:sandwich}
 \Psi(v,\mu)+c\norm y_{H^1}^2
 \le\mathcal E(w,\mu)\le
 \Psi(v,\mu)+C\norm y_{H^1}^2.
\end{equation}
A positive definite $Q_\rho$ then makes the critical ball
a strict local minimizer among small $C^1$ equal-volume graphs.
\end{lemma}
\begin{proof}
The exact area and interaction first variations are
\begin{equation}\label{eq:area-gradient}
 P(E_w)=\int r_w^2\sqrt{1+r_w^{-6}|\nabla w|^2},\qquad
 DI_\rho(w)[\eta]=\int\rho(r_w)u_w(r_w\omega)\eta,\quad
 u_w=G_\kappa*(\rho\one_{E_w}).
\end{equation}
The projected Euler operator is smooth from mean-zero $C^{2,\gamma}$
to mean-zero $C^{0,\gamma}$. For the nonlocal term, we account for the moving singularity as follows. Extend the boundary map to a
smoothly $w$-dependent, uniformly bi-Lipschitz map $\Phi_w$ of space,
equal to the identity near zero and outside a fixed annulus. The
potential trace becomes
\[
 \int_B G_\kappa(\Phi_w(\omega)-\Phi_w(y))
 \rho(\Phi_w(y))\det D\Phi_w(y)\dd y.
\]
Each parameter derivative of the kernel contains differences of
parameter derivatives of $\Phi_w$.
Each such difference is $O(|\omega-y|)$, with a bound multilinear in
the $C^{2,\gamma}$ norms of the directions. So the differentiated
integrands are bounded by $C|\omega-y|^{-1}$; one tangential derivative
is bounded by $C|\omega-y|^{-2}$. Both are uniformly integrable in
three dimensions. Taylor's formula along parameter segments gives the
same bounds for Fr\'echet remainders, with a factor tending to zero.
Integration then controls these remainders in the Lipschitz norm
of the trace, and so in $C^{0,\gamma}$. This proves smoothness to
every finite order. Smoothness of the local area operator is immediate
from \eqref{eq:area-gradient}.

The linearization on $W$ is $-\Delta_{\Sph}-2$ plus a constant
and a lower-order single-layer operator. It is Fredholm of index zero
between these H\"older spaces, and the spectral gap makes it invertible
by elliptic estimates. The implicit function theorem gives $z$.
The ball is critical for all $\mu$ and the linearization preserves
$V$ and $W$, proving the two vanishing identities. Uniqueness gives
equivariance. Central inversion acts as $-1$ on $V_3$, so $\Psi$ is
even. Its quadratic term is exactly $-\alpha\mu\norm v_2^2/2$,
which proves \eqref{eq:reduced} by Taylor expansion.

To compare $C^1$ graphs, we prove that the Hessian is continuous as an
$H^1$ bilinear form.
Differentiating \eqref{eq:area-gradient} gives, with $X_w=r_w\omega$,
\begin{align}\label{eq:second-interaction}
 D^2I_\rho(w)[\xi,\eta]
 ={}&\int\frac{\rho'(r_w)u_w(X_w)+
 \rho(r_w)\omega\cdot\nabla u_w(X_w)}{r_w^2}\xi\eta\dd\omega\notag\\
 &+\iint\rho(r_w(\omega))\rho(r_w(\sigma))
 G_\kappa(X_w(\omega)-X_w(\sigma))\xi(\omega)\eta(\sigma)
 \dd\omega\dd\sigma.
\end{align}
The derivative in the first line is radial. The volume potential has
a continuous gradient across the boundary. If $\delta=\norm w_{C^0}$,
then the source difference is bounded and supported on a set of
volume $O(\delta)$; the kernels $|x|^{-1}$ and $|x|^{-2}$ give
$\norm{u_w-u_0}_\infty\le C\delta^{2/3}$ and
$\norm{\nabla u_w-\nabla u_0}_\infty\le C\delta^{1/3}$.
The first coefficient in \eqref{eq:second-interaction} then
converges uniformly. For the second line,
\[
 |X_w(\omega)-X_w(\sigma)|\ge\min r_w\,|\omega-\sigma|,
 \quad
 |(X_w-\mathrm{Id})(\omega)-(X_w-\mathrm{Id})(\sigma)|
 \le C\norm w_{C^1}|\omega-\sigma|.
\]
The difference from the ball kernel is at most
$C\norm w_{C^1}/|\omega-\sigma|$. Its $L^2$ operator norm tends
to zero by Schur's bound. The perimeter Hessian has coefficients
continuous in $w,\nabla w$. Combining these estimates,
\begin{equation}\label{eq:two-norm}
 |(D^2\mathcal E(w,\mu)-D^2\mathcal E(0,0))[\xi,\eta]|
 \le o(1)\norm\xi_{H^1}\norm\eta_{H^1}
\end{equation}
as $\norm w_{C^1}+|\mu|\to0$. The integral formulas and Taylor's
formula extend to $C^1$ graphs by smooth approximation and the displayed
uniform bounds.

The spectral gap is coercive on $W$ in $H^1$, so
\eqref{eq:two-norm} preserves that coercivity along the segment
$v+z(v,\mu)+sy$. The first derivative along $y$ vanishes at $s=0$.
Integrating the second derivative proves \eqref{eq:sandwich}.
If $Q_\rho(v)\ge a\norm v_2^4$, then at $\mu=0$ this gives
$\mathcal E(w,0)\ge(a/2)\norm v_2^4+c\norm y_{H^1}^2$ after
shrinking the neighborhood. Equality forces $w=0$.
\end{proof}

The graph topology and the coercivity norm differ, as in the shape
stability results of \cite{DL}. Here the needed continuity estimate is
\eqref{eq:two-norm}.

\subsection{One positive radial coefficient and one angular control}
Set $\N(v)=\norm v_2^4$ and $\J(v)=\int_{\Sph}v^4$. The invariant
calculation in Section~\ref{sec:selection} shows that every quartic
$O(3)$ invariant on $V_3$ is uniquely $A\N+B\J$.

\begin{proposition}[The two quartic forms needed for the construction]
\label{prop:quartic}
For every sufficiently small fixed positive $\kappa$ there exist two
densities satisfying Lemma~\ref{lem:jets} whose reduced quartics are
\begin{equation}\label{eq:two-quartics}
 Q_T=A_\kappa\N+b_\kappa\J,\qquad
 Q_S=A_\kappa\N-b_\kappa\J,
 \quad A_\kappa>0,\quad 0<b_\kappa<A_\kappa/j_S,
 \qquad j_S=\frac{175}{286\pi}.
\end{equation}
Both can be realized while fixing
$\rho''(1)=\rho_0''(1)+1$.
\end{proposition}
\begin{proof}
First fix this second derivative and set $\rho'''(1)=0$, using
Lemma~\ref{lem:jets}. Finite-order shape derivatives depend only on
the boundary derivatives of the required orders and the ball Yukawa
moment. This also follows from the signed-shell argument below.
The coefficient $A_\kappa$ of $\N$ does not depend on the
particular interior correction. The direct transmission calculation
in Appendix~\ref{app:calculation} gives
\begin{equation}\label{eq:Apositive}
 \lim_{\kappa\downarrow0}A_\kappa
 =\frac{66345100\sqrt7-84268315}{1672704\pi}>0.
\end{equation}
It is positive because $\sqrt7>2$. The coefficients
involve only finitely many homogeneous radial solutions and finitely
many boundary derivatives, and are continuous as $\kappa\downarrow0$.
So $A_\kappa>0$ for small positive $\kappa$.

To change the third derivative by $\eta$, use another density with
the same second and lower derivatives and the same moment. Write its
difference as $d$, let $h=r_w-1$, $\sigma_h=\one_{E_w}-\one_B$,
and put $u_d=G_\kappa*(d\one_B)$. Since $M_\kappa(d)=0$,
$u_d(1)=u_d'(1)=0$. The radial equation
$(ru_d)''-\kappa^2(ru_d)=-rd$ implies
$u_d(1+s)=O(|s|^5)$ inside, while it vanishes outside. The exact
difference of the relative interaction energies is
\begin{equation}\label{eq:shell}
 \int\sigma_h\{d u_B+(\rho+d)u_d\}
 +\ip{d\sigma_h}{G_\kappa*(\rho\sigma_h)}
 +\frac12\ip{d\sigma_h}{G_\kappa*(d\sigma_h)}.
\end{equation}
Here $d(1+s)=\eta s^3/6+O(s^4)$. The first contribution is
$U\eta\int h^4/24+O(\norm h_\infty^5)$. A spherical shell of
thickness $\delta$ has Coulomb potential $O(\delta)$ uniformly;
since $G_\kappa\le G_0$, the cross terms in \eqref{eq:shell} are
$O(\norm h_\infty^5)$. So the first changed energy term is
$\beta_*U\eta\J(w)/24$ at degree four. More generally the same
argument shows that if the boundary derivatives through order three
and the moment agree, all derivatives of relative energy through
order four agree.

The quadratic complement correction is determined by the cubic
energy, which has not changed. The two reduction graphs
first differ at order three, and the energy cost of that difference
is order six. The change in the reduced quartic is
\begin{equation}\label{eq:third-jet}
 Q_{\rho+d}-Q_\rho=\frac{\beta_*U}{24}\eta\J.
\end{equation}
Since $\beta_*U>0$, this sets its $\J$ coefficient to either
$b_\kappa$ or $-b_\kappa$, for example with
$b_\kappa=A_\kappa/(2j_S)$. Both prescribed third derivatives are
finite for the fixed $\kappa$ and are realized by Lemma~\ref{lem:jets}.
We take the limit only in the local coefficient formula. The interior
corrections need not depend continuously on $\kappa$ at zero.
\end{proof}

\section{Angular selection and the unique near-ball minimizing orbit}
\label{sec:selection}

We compute the angular bounds and Hessians for harmonic cubics in our
normalization. The invariant theory and orbit types are described in
\cite{IG,CLM,CHQZ,BvV}.

\subsection{Two invariants and their extremal orbits}
Write $v=T_{ijk}x_ix_jx_k$ with $T$ symmetric and trace-free, and set
\[
 A_{ij}=T_{ik\ell}T_{jk\ell},\qquad I=(\tr A)^2,\qquad K=\tr(A^2),
 \qquad Q_\ell(v)=\norm{\Pi_\ell(v^2)}_2^2.
\]
Repeated indices are summed.

\begin{lemma}[Quartic invariants, equality orbits, and angular Hessians]
\label{lem:angular}
The quartic $O(3)$ invariants on $V_3$ have basis $(\N,\J)$, and
\begin{equation}\label{eq:angular-range}
 j_T:=\frac{315}{572\pi}\le\J(v)\le
 \frac{175}{286\pi}=:j_S\qquad(\norm v_2=1).
\end{equation}
The equality sets are respectively $\SO e_T$ and $\SO e_S$.
On the three-dimensional normal space to either orbit inside the
unit sphere of $V_3$, the Hessian of $\J$ has eigenvalues
\begin{equation}\label{eq:angular-hessian}
 e_T:\quad\frac{56}{143\pi}\ (\text{triple}),\qquad
 e_S:\quad-\frac{63}{143\pi}\ (\text{simple}),\quad
 -\frac{42}{143\pi}\ (\text{double}).
\end{equation}
\end{lemma}
\begin{proof}
The complexified $V_3$ has rotation weights $-3,-2,\ldots,3$.
In its fourth symmetric power the multiplicities at weights zero
and one are respectively $18$ and $16$: count nondecreasing
quadruples in $\{0,\ldots,6\}$ of sums $12$ and $13$.
Their difference, two, is the multiplicity of the trivial representation. Quartics are even under central inversion, so the
$SO(3)$ and $O(3)$ invariant spaces coincide. The invariants $I,K$
are independent, having ratios $K/I=1/3$ for $xyz$ and $1/2$ for
$z^3-3zy^2$.

The moment formula
\begin{equation}\label{eq:moments}
 \int_{\Sph}x^{2a}y^{2b}z^{2c}\dd\omega
 =4\pi\frac{(2a-1)!!(2b-1)!!(2c-1)!!}
 {(2a+2b+2c+1)!!}
\end{equation}
and vanishing of odd-coordinate moments give
\begin{equation}\label{eq:tensor-integrals}
 \norm v_2^2=\frac{8\pi}{35}\tr A,\quad
 \J(v)=\frac{16\pi}{5005}(7I+6K),\quad
 Q_2=\frac{32\pi}{735}(K-I/3),\quad
 \J=j_T\N+\frac{63}{143}Q_2.
\end{equation}
For example $\Pi_2(v^2)=\frac27(A-\frac13\tr A\,\mathrm{Id})_{ij}x_ix_j$;
this follows by applying the Laplacian twice to the homogeneous
polynomial $v^2$. These identities prove that $(\N,\J)$ is a basis.

The matrix $A$ is nonnegative, so $K\ge I/3$. For the upper bound
choose a positive maximum of $v$ as north pole and diagonalize its
tangential quadratic part. Then $T_{133}=T_{233}=T_{123}=0$.
Put $a=T_{333}>0$, $b=T_{113}$, $p=T_{111}$, $q=T_{112}$,
$s=p^2+q^2$. The trace conditions give
\[
 A=\begin{pmatrix}
 2(s+b^2)&0&p(a+2b)\\
 0&2(s+(a+b)^2)&q(a+2b)\\
 p(a+2b)&q(a+2b)&a^2+b^2+(a+b)^2
 \end{pmatrix},\qquad
 I-2K=12\{b^2(a+b)^2+a^2(p^2+q^2)\}.
\]
So $K\le I/2$, with equality exactly when $p=q=0$ and
$b=0$ or $b=-a$; both give rotations and multiples of
$z^3-3zy^2$. If $A$ is scalar, then $b=-a/2$ and
$p^2+q^2=a^2/2$. A rotation about the pole makes
$q=0$, $p=a/\sqrt2$, giving
\[
 a\left[z^3-\frac32(x^2+y^2)z+\frac1{\sqrt2}(x^3-3xy^2)\right]
 =3\sqrt3\,aXYZ,
\]
where $X=\sqrt{2/3}x+z/\sqrt3$,
$Y=-x/\sqrt6+y/\sqrt2+z/\sqrt3$, and
$Z=-x/\sqrt6-y/\sqrt2+z/\sqrt3$ are orthogonal coordinates.
This proves the lower equality classification. Formula
\eqref{eq:tensor-integrals} now gives \eqref{eq:angular-range}.
The negatives of both representatives are proper rotations of them:
combine a coordinate interchange with central inversion for $e_T$,
and use a $\pi/3$ rotation about $z$ for $e_S$.

At a unit critical harmonic $e$ the constrained angular Hessian is
\begin{equation}\label{eq:hessian-moment}
 H_e[h,h]=12\int e^2h^2-4\J(e)\norm h_2^2,
 \qquad\ip e h=0.
\end{equation}
At $e_T$, a transverse orthonormal basis consists of
$\sqrt7\,x(5x^2-3|x|^2)/(4\sqrt\pi)$ and its cyclic partners.
At $e_S$, use
\[
 \frac{\sqrt7}{4\sqrt\pi}(5z^3-3z|x|^2),\qquad
 \sqrt{\frac{21}{32\pi}}x(5z^2-|x|^2),\qquad
 \sqrt{\frac{21}{32\pi}}y(5z^2-|x|^2).
\]
The corresponding diagonal moments $\int e^2h^2$ are
\[
 e_T:\quad\frac{371}{1716\pi}\ \text{in each direction};\qquad
 e_S:\quad\frac{287}{1716\pi},\ \frac7{39\pi},\ \frac7{39\pi}.
\]
All mixed moments vanish by parity. Equation
\eqref{eq:hessian-moment} gives \eqref{eq:angular-hessian}.
The three rotational directions are exact zero modes by invariance.
Their polynomial spans are $x(z^2-y^2),y(x^2-z^2),z(y^2-x^2)$ at
$e_T$, and $3x^2y-y^3,xyz,z(x^2-y^2)$ at $e_S$. These are
orthogonal to the transverse vectors displayed above.
Since the unit sphere has dimension six and each orbit dimension
three, no further directions remain.
\end{proof}

\subsection{From quartic preference to exact minimizing shapes}

\begin{proposition}[Unique near-ball minimum modulo rotations]
\label{prop:selection}
Suppose the reduced quartic is one of \eqref{eq:two-quartics}, and
let $e=e_T$ or $e_S$ for the chosen quartic. Put
$a=Q(e)>0$. At $\mu=0$ the ball is a strict $C^1$ local minimizer.
For every sufficiently small $\mu>0$, the shapes attaining the
least energy in a fixed sufficiently small $C^1$ graph neighborhood
of the ball form exactly one rotation orbit, with representative
\[
 w=t e+z(te,\mu),\qquad
 t^2=\frac{\alpha\mu}{4a}+O(\mu^2).
\]
Its energy is $-\alpha^2\mu^2/(16a)+O(\mu^3)$ relative to the ball.
Its full stabilizer is that of $e$, and it is strictly convex for
small $\mu$.
\end{proposition}
\begin{proof}
Lemma~\ref{lem:angular} shows that $Q(v)\ge a\norm v_2^4$ with
equality precisely along the chosen angular orbit. Positivity at
$\mu=0$ follows from Lemma~\ref{lem:reduction}.
The fixed-point space of the stabilizer of $e$ in $V_3$ is $\R e$.
For $T_d$, the even coordinate sign changes and permutations leave
only $xyz$; for $D_{3h}$, reflection in $z$ and the planar rotation
and reflection leave only the cosine mode of frequency three.
By equivariance, a stationary point on that line is a
stationary point of the full reduced energy. If
$\psi(t,\mu)=\Psi(te,\mu)$, its equation is
\[
 \frac{\partial_t\psi(t,\mu)}t
 =-\alpha\mu+4at^2+O(t^4+|\mu|t^2)=0.
\]
The implicit function theorem in $\mu$ gives
$\mu(t)=4at^2/\alpha+O(t^4)$. This function is even and strictly
increasing for small $t>0$, so every sufficiently small $\mu>0$
has exactly one small positive amplitude, with the stated expansion.

To prove angular uniqueness near the selected orbit, take a normal slice
to the orbit in the finite-dimensional unit sphere of $V_3$. Divide the
angular equation at $v=t\omega$ by $t^4$. It extends smoothly to
$t=\mu=0$, with invertible angular derivative given by
$\pm b_\kappa H_e$. Symmetry makes the slice origin a solution for
all small $t,\mu$. The implicit function theorem makes it the unique
one there. The radial equation above then gives a single critical
orbit in the cone. The slice argument takes place in $V_3$, so no differentiability of the
rotation action on the H\"older graph space is needed.

We now show that every near-ball minimizer lies in that cone.
Choose $\delta>0$ and then $\mu_0>0$ so that the remainder in
\eqref{eq:reduced} is at most $(a/2)\norm v_2^4$ for
$\norm v_2\le\delta$, $0<\mu<\mu_0$. Then
\[
 \Psi(v,\mu)\ge-\frac\alpha2\mu\norm v_2^2+
 \frac a2\norm v_2^4.
\]
Choose $K^2>\alpha/a$ and decrease $\mu_0$ so that
$K\sqrt{\mu_0}<\delta$. The energy is positive when
$K\sqrt\mu\le\norm v_2\le\delta$, while the constructed branch
has negative energy. The minimum on the closed $\delta$ ball is
then interior and has $\norm v_2<K\sqrt\mu$.
Under $v=\sqrt\mu\,y$, the energies divided by $\mu^2$ converge
in $C^2$ on that fixed $y$ ball to
\[
 \Phi(y)=-\frac\alpha2\norm y_2^2+Q(y).
\]
Its global minimizing set is exactly
$\sqrt{\alpha/(4a)}\,\SO e$. Any sequence of minimizers of the
scaled energies has a convergent subsequence; comparison with each
fixed $y$ makes its limit a minimum of $\Phi$. It follows, uniformly
over all minimizers, that they enter the previously chosen cone and
stay away from zero in scaled coordinates. Angular and radial
uniqueness identify them with the constructed orbit. This step uses the angular bounds on the whole unit sphere.

Choose the fixed $C^1$ graph neighborhood so that its center projections
lie in the $\delta$ ball and \eqref{eq:sandwich} holds. That estimate
compares every graph in the neighborhood with its exact reduction graph.
Each graph has energy at least the reduced minimum, with equality only
if $y=0$ and its $V_3$ coordinate
lies on that orbit. The constructed graphs lie in the
chosen neighborhood for small $\mu$, so they attain the minimum.
Substitution in \eqref{eq:reduced} gives the stated energy expansion.
Equivariance gives every symmetry of $e$ as a symmetry of the graph;
projecting the graph onto $V_3$ shows that any symmetry of the graph
also fixes $e$. Finally $z(te,\mu)=O_{C^{2,\gamma}}(t^2)$ and $t\to0$,
so the boundaries are nonspherical strictly convex perturbations of
the unit sphere.
\end{proof}

\section{From local selection to the first global transition}\label{sec:global}

Fix $\kappa>0$ and write $\mathsf V=|B|$. In this section we consider
smooth positive densities satisfying
\begin{equation}\label{eq:quadraticgrowth}
 c_\rho(1+|x|^2)\le\rho(x)\le C_\rho(1+|x|^2),
 \qquad x\in\R^3,
\end{equation}
with positive constants depending on the density. All densities constructed
above satisfy these bounds. So do their modifications by any fixed bounded
nonnegative function. Sets are identified up to null sets.
To pass from graph minimality to $L^1$ minimality, we use penalization
as in \cite{AFM,BC}. We give the argument for a Hessian with a kernel.
For BV compactness and perimeter slicing, we use \cite{Maggi}.

\begin{lemma}[Existence and parameter compactness]\label{lem:compact}
For every $\beta>0$, $F_{\beta,\rho}$ has a global minimizer among
finite-perimeter sets of volume $\mathsf V$.
If $\beta_j\to\beta_0>0$ and $E_j$ are corresponding global minimizers,
a subsequence converges in $L^1(\R^3)$ to a global minimizer at $\beta_0$.
\end{lemma}
\begin{proof}
For $z\in\mathbb Z^3$, let $Q_z=z+[-1/2,1/2)^3$,
$m_z=\inf_{Q_z}\rho$, and
$g=e^{-\kappa\sqrt3}/(4\pi\sqrt3)$.
Positivity of the kernel gives
\[
 2I_\rho(E)\ge g\sum_zm_z^2|E\cap Q_z|^2,
 \qquad
 \sum_{z\in\mathcal Z}|E\cap Q_z|
 \le\left(\frac{2I_\rho(E)}g\right)^{1/2}
      \left(\sum_{z\in\mathcal Z}m_z^{-2}\right)^{1/2}.
\]
The lower bound in \eqref{eq:quadraticgrowth} implies
$\sum_zm_z^{-2}<\infty$. Bounded self-energy then gives uniform
tightness of the ordinary volume. Comparison with $B$ bounds the perimeter
and self-energy of a minimizing sequence. Local BV compactness and tightness
give global $L^1$ convergence preserving volume. Perimeter is lower
semicontinuous, and Fatou's lemma treats the positive double integral.
This proves existence.

For varying parameters, comparison with $B$ and $P(E_j)\ge P(B)$ give
$I_\rho(E_j)\le I_\rho(B)$ and a uniform perimeter bound. The same compactness argument applies. Since $I_\rho(E_j)$ is bounded,
replacing $\beta_j$ by $\beta_0$ changes the energy by $o(1)$. Lower
semicontinuity and comparison with each fixed finite-energy competitor
then prove that the limit is a minimizer. Infinite-energy
competitors impose no additional condition at a positive coupling.
\end{proof}

We need the following local estimate for the unbounded density. For a compact set $K$ and any $|E|\le\mathsf V+1$,
\begin{equation}\label{eq:localpotential}
 \sup_{x\in K}G_\kappa*(\rho\one_E)(x)\le C_{K,\mathsf V}.
\end{equation}
The contribution from $|x-y|<1$ is bounded by local boundedness of
$\rho$ and integrability of $G_\kappa$. On $|x-y|\ge1$, the function
$\rho(y)G_\kappa(x-y)$ is uniformly bounded for $x\in K$, by the upper
quadratic bound and exponential decay. Its integral over $E$ is bounded by
that supremum times $|E|$.
So if $I_\rho(E)<\infty$ and $F\triangle E\subset K$, then
\begin{equation}\label{eq:localchange}
 |I_\rho(F)-I_\rho(E)|\le C_{K,\mathsf V}|F\triangle E|.
\end{equation}
To see this, set $a=\rho(\one_F-\one_E)$ and expand the difference as
$\int a(G_\kappa*(\rho\one_E))+\frac12\int a(G_\kappa*a)$.
The first term is controlled by \eqref{eq:localpotential}; the second by
$\frac12(\sup_K\rho)^2\|G_\kappa\|_1|F\triangle E|$.
So bounded local modifications preserve finiteness of the energy.

\begin{lemma}[Uniform regularity near the ball]\label{lem:regularity}
Let $\beta_j$ lie in a compact subinterval of $(0,\infty)$, and let $H_j$
be real-valued penalties with a common Lipschitz constant for the distance
$|E\triangle F|$.
Suppose $A_j$ minimize $F_{\beta_j,\rho}+H_j$ at volume $\mathsf V$ and
$A_j\to B$ in $L^1$.
For all sufficiently large $j$, their essential supports are uniformly
bounded, and there exist $\Lambda,r_0>0$ such that
\begin{equation}\label{eq:areaminimality}
 P(A_j;B_r(x))\le P(F;B_r(x))+\Lambda|F\triangle A_j|
\end{equation}
whenever $r\le r_0$ and $F\triangle A_j\Subset B_r(x)$.
Their boundaries are eventually spherical graphs converging to $\partial B$
in $C^{1,\alpha}$ for every $0<\alpha<1/2$.
\end{lemma}
\begin{proof}
Comparison with $B$ bounds $P(A_j)$ and $I_\rho(A_j)$ uniformly, since
$|H_j(A_j)-H_j(B)|\le2\mathsf V\operatorname{Lip}(H_j)$.
Choose bounded neighborhoods $U_1,U_2$ with disjoint closures near two
separated sphere patches, and smooth fields $X_i$ compactly supported in
$U_i$, with $b_i=\int_B\operatorname{div}X_i>0$.
Then $\int_{A_j}\operatorname{div}X_i\ge b_i/2$ for large $j$.
For the flows $\Phi_{i,t}$, bounded second derivatives of the Jacobian
give a common interval on which the derivative of
$|\Phi_{i,t}(A_j)|$ stays positive. Any volume correction $|s|<s_0$
can be made with $|t|\le C|s|$.
The tangential-Jacobian and local BV flow estimates give
\[
 |P(\Phi_{i,t}(A_j))-P(A_j)|
   +|\Phi_{i,t}(A_j)\triangle A_j|
 \le C|t|P(A_j;U_i).
\]
Equations \eqref{eq:localpotential}--\eqref{eq:localchange} and the penalty
bound show that the full correction costs at most $C|s|$.
If a preliminary modification leaves $U_i$ unchanged, the volume-change
function and local perimeter in this argument are unchanged. The same
cost estimate applies to the truncations and bounded local modifications
below; their intermediate volumes are at most $\mathsf V+1$.

For $t$ outside the two correction neighborhoods, put
$m_j(t)=|A_j\setminus B_t|$ and
$a_j(t)=-m'_j(t)=\mathcal H^2(A_j\cap\partial B_t)$ for almost every $t$.
The perimeter slicing identity reads
\[
 P(A_j\cap B_t)-P(A_j)=2a_j(t)-P(A_j\setminus B_t).
\]
Removing the tail decreases the positive interaction energy. Whenever
$m_j(t)<s_0$, restore its volume using either flow. Minimality, the correction estimate, and the
isoperimetric inequality yield
\[
 c_3m_j(t)^{2/3}\le P(A_j\setminus B_t)
                 \le2a_j(t)+Cm_j(t).
\]
Choose a fixed $T$ containing $B,U_1,U_2$. Since $m_j(T)\to0$, the correction
is allowed and $Cm_j\le(c_3/2)m_j^{2/3}$ for $t\ge T$ and large $j$.
Where $m_j>0$, it follows that
$-(m_j^{1/3})'\ge c_3/12$. So all these tails vanish before a common
finite radius $L$.

Choose $r_0$ smaller than one quarter of the distance between
$\overline U_1,\overline U_2$, with $|B_{r_0}|<\min(s_0,1)$.
If $B_r(x)$ is disjoint from $B_L$, \eqref{eq:areaminimality} is automatic.
Otherwise all changes occur in a fixed enlarged ball. One $U_i$ is disjoint
from $B_r(x)$; use its flow to correct the volume error of $F$, whose size
is at most $|F\triangle A_j|$.
The local perimeter on $U_i$ is still that of $A_j$, even when $P(F)$ is
large. The correction costs and \eqref{eq:localchange} then give
$P(A_j)\le P(F)+C|F\triangle A_j|$ by minimality.
The exterior perimeter terms cancel, proving \eqref{eq:areaminimality}.

Since $|F\triangle A_j|\le|B_r|$, this is area almost minimality with a
common error $C r^3$. We apply the spherical convergence theorem in
\cite[Theorem~2.2]{FJ}: uniform perimeter bounds and global
$L^1$ convergence to $B$ imply eventual $C^{1,1/2}$ spherical graphs and
convergence of their graph functions in $C^{1,\alpha}$, $\alpha<1/2$.
\end{proof}

\begin{lemma}[Strict graph minimality implies strict $L^1$ minimality]
\label{lem:L1}
Fix $\beta>0$. If $B$ is a strict minimizer among sufficiently small
$C^1$ normal graphs of volume $\mathsf V$, it is a strict $L^1$ local minimizer
among all finite-perimeter sets of that volume. The constrained Hessian may have a kernel.
\end{lemma}
\begin{proof}
Choose a smooth compactly supported field $X$ with $|X|\le1$ and
$X=\nu_B$ on $\partial B$. The divergence formula gives
$P(E)-P(B)\ge-C|E\triangle B|$.
Positivity also gives $I_\rho(E)\ge I_\rho(E\cap B)$, while
$I_\rho(B)-I_\rho(E\cap B)\le C|B\setminus E|$.
Writing $d(E)=|E\triangle B|$, we get
\[
 F_{\beta,\rho}(E)\ge F_{\beta,\rho}(B)-Cd(E).
\]
If the assertion fails, choose $0<\delta_j=d(E_j)\to0$ with
$F_{\beta,\rho}(E_j)\le F_{\beta,\rho}(B)$.
For a fixed $M>C$, minimize
$J_j(E)=F_{\beta,\rho}(E)+M|d(E)-\delta_j|$ at volume $\mathsf V$.
The direct method above applies to this nonnegative $L^1$-continuous
penalty. Its minimizers satisfy
\[
 J_j(A_j)\le F_{\beta,\rho}(B)<J_j(B),\qquad
 d(A_j)\le\frac{M\delta_j}{M-C},\qquad
 F_{\beta,\rho}(A_j)\le F_{\beta,\rho}(B).
\]
The penalties have common Lipschitz constant $M$.
Lemma~\ref{lem:regularity} makes the nonball sets $A_j$ eventually small
$C^1$ graphs of the prescribed volume, a contradiction.
\end{proof}

\begin{proposition}[A common finite exterior density barrier]\label{prop:wall}
Consider finitely many radial densities satisfying \eqref{eq:quadraticgrowth}.
Suppose that, at a common $\beta_*>0$, the ball is a strict $C^1$ graph
minimizer in each model, while every slightly larger coupling admits
arbitrarily small equal-volume graphs of lower energy than the ball.
There exist $c>1$, a bounded smooth nondecreasing radial function $h_c$,
zero on $B_c$ and positive outside, and a common finite $s\ge0$, such that
the modified densities $\widehat\rho_i=\rho_i+s h_c$ have
$\beta_g(\widehat\rho_i)=\beta_*$.
The centered ball is their unique global minimizer for
$0<\beta\le\beta_*$.
All global minimizers converge to $B$ in $L^1$ as $\beta\downarrow\beta_*$,
uniformly over the choice of minimizer.
The barrier preserves the exact energy near the ball, all ball data, and
strict radial monotonicity when present.
\end{proposition}
\begin{proof}
Lemma~\ref{lem:L1} provides strict $L^1$ radii; let $\delta$ be their minimum.
Choose $c>1$ with $2(|B_c|-\mathsf V)<\delta$.
Every $E\subset B_c$ of volume $\mathsf V$ then satisfies $|E\triangle B|<\delta$,
so $B$ is the unique minimizer in this hard-wall class.
Fix $h_c$, for example $h_c(r)=e^{-(r-c)^{-2}}$ for $r>c$ and zero otherwise.

Work first with one density. If arbitrarily large finite strengths fail,
there are $s_j\to\infty$ and nonball global minimizers $E_j$ for
$\rho_j=\rho+s_jh_c$ at $\beta_*$.
Existence holds for each fixed $s_j$. Comparison with $B$ and the
isoperimetric inequality give a uniform perimeter bound and
\[
 I_\rho(E_j)\le I_{\rho_j}(E_j)\le I_\rho(B),\qquad
 s_j^2 I_{h_c}(E_j)\le I_\rho(B).
\]
The tightness estimate for the original density gives global $L^1$
convergence along a subsequence. Fatou's lemma forces its limit $E_\infty$
to satisfy $I_{h_c}(E_\infty)=0$, so $E_\infty\subset B_c$ by strict
positivity of the kernel. Lower semicontinuity gives
$F_{\beta_*,\rho}(E_\infty)\le F_{\beta_*,\rho}(B)$, so $E_\infty=B$.
Strict $L^1$ minimality for the original density now contradicts
\[
 F_{\beta_*,\rho_j}(E_j)\ge F_{\beta_*,\rho}(E_j)
 >F_{\beta_*,\rho}(B)=F_{\beta_*,\rho_j}(B).
\]
This proves that every sufficiently large finite strength works. Take the largest
required threshold for the finitely many models and fix one common $s$.
The contradiction uses strict $L^1$ minimality for the fixed original
density, so it needs no regularity bound uniform in $s$.

For a finite-energy competitor and $0<\beta<\beta_*$,
\[
 \Delta F_\beta
 =\frac\beta{\beta_*}\Delta F_{\beta_*}
  +\left(1-\frac\beta{\beta_*}\right)\Delta P
\]
gives strict uniqueness, since $\Delta P\ge0$.
Infinite-energy competitors are irrelevant for positive couplings.
Graphs sufficiently close to $B$ lie in $B_c$, so their exact energies
are unchanged. The lower-energy competitors just above $\beta_*$ persist.
For each such competitor at $\beta_1>\beta_*$,
$\Delta P+\beta_1\Delta I<0$ and $\Delta P\ge0$ imply $\Delta I<0$;
the same competitor defeats the ball at every $\beta\ge\beta_1$.
So the ball cannot become minimizing again, and
$\beta_g=\beta_*$.
Parameter compactness and uniqueness at $\beta_*$ give the stated
uniform convergence. Finally, $h_c=0$ on $B_c$ preserves local energies,
and $h_c'\ge0$ preserves strict radial monotonicity.
The common addition also preserves agreement of the densities at all
sufficiently large radii. Its fixed finite size preserves
\eqref{eq:quadraticgrowth}.
\end{proof}

\begin{proposition}[Every global minimizer lies on the exact reduction graph]
\label{prop:globalgraph}
Fix one model after the finite barrier construction, with the reduction
and $C^1$ energy estimate of Lemma~\ref{lem:reduction}.
For all sufficiently small $\mu>0$, every global minimizer has the form
$E_{v+z(v,\mu)}$, where $v$ minimizes $\Psi(\cdot,\mu)$ over a fixed small
center neighborhood.
\end{proposition}
\begin{proof}
For any sequence $\mu_j\downarrow0$, global minimizers converge to $B$
in $L^1$ by Proposition~\ref{prop:wall}.
Lemma~\ref{lem:regularity}, with zero penalty, gives small $C^1$ graphs.
In exact volume coordinates write $w_j=v_j+z(v_j,\mu_j)+y_j$, $y_j\in W$.
Compare the minimizer with the reduced shape at the same $v_j$.
Equation~\eqref{eq:sandwich} gives
\[
 \Psi(v_j,\mu_j)+c\|y_j\|_{H^1}^2
 \le F(E_j)-F(B)\le\Psi(v_j,\mu_j).
\]
The estimate forces $y_j=0$. Comparison with every other small reduced shape proves
minimality in the center neighborhood. A sequence of counterexamples
would contradict this argument, so the conclusion holds for every
minimizer throughout a sufficiently small parameter interval.
\end{proof}

If the reduced minimizers form one rotation orbit, existence and rotation
invariance imply that the global minimizers form the same orbit. All
constants in this section concern fixed densities. We first fix the
boundary derivatives, then the finite exterior barrier, and then shrink
the parameter interval.

\section{Proof of the main theorem}\label{sec:completion}

\begin{proof}[Proof of Theorem~\ref{thm:main}]
Choose $\kappa_0$ small enough for Lemma~\ref{lem:crossing} and
Proposition~\ref{prop:quartic}. Fix $0<\kappa<\kappa_0$ for the
rest of the construction. Use Lemma~\ref{lem:jets} to realize
\eqref{eq:two-quartics}, with the same second derivative
$\rho''(1)=4-\sqrt7$ and the two prescribed third derivatives.
Call these preliminary densities $\rho_T^0,\rho_S^0$. They preserve
the common boundary value, first derivative, mass, Yukawa moment,
and ball self-energy. Lemma~\ref{lem:spectrum} gives identical
constrained Hessians for every coupling, including all degree-one
and all high-degree modes. Both densities satisfy two-sided
quadratic growth bounds and are strictly increasing.

Their quartics have positive angular minima
\[
 a_T=A_\kappa+b_\kappa j_T>0,\qquad
 a_S=A_\kappa-b_\kappa j_S>0.
\]
Lemma~\ref{lem:reduction} gives strict $C^1$ graph minimality of the
critical ball, and Lemma~\ref{lem:L1} upgrades it to strict $L^1$
local minimality. Choose the common exterior radius and a common
finite barrier strength as in Proposition~\ref{prop:wall}, and set
\[
 \rho_T=\rho_T^0+s h_c,\qquad
 \rho_S=\rho_S^0+s h_c.
\]
The increment is smooth, nonnegative, nondecreasing, bounded, and
zero throughout $B_c$, where $c>1$. It preserves strict increase,
quadratic growth, all the ball data, and the entire exact energy of
every shape contained in $B_c$. The two densities coincide away
from the compact annuli supporting their preliminary perturbations.

Proposition~\ref{prop:wall} makes $B$ uniquely globally minimizing
for $0<\beta\le\beta_*$. Proposition~\ref{prop:selection} gives
near-ball competitors of lower energy for every sufficiently small
$\mu>0$; these lie in $B_c$ and remain unchanged by the barrier.
Proposition~\ref{prop:wall} then gives $\beta_g=\beta_*$. The common linear threshold and its limit
follow from Lemma~\ref{lem:crossing}.

Fix these two final densities. By compactness and uniqueness at
$\beta_*$, all their global minimizers converge to $B$ in $L^1$
as $\beta\downarrow\beta_*$. Proposition~\ref{prop:globalgraph}
places every such minimizer on the exact reduction graph. Its
center coordinate minimizes the reduced energy in a fixed small
neighborhood. Proposition~\ref{prop:selection} identifies that
minimizing set as the single prescribed rotation orbit. Global
existence gives a point of each orbit, and rotation invariance
makes the whole orbit globally minimizing. Taking the smaller of
the two parameter neighborhoods gives the common $\varepsilon$.
The stabilizers, convexity, and expansions
\eqref{eq:main-asymptotic} follow from the same proposition.
\end{proof}

\appendix

\section{A finite transmission calculation for the positive coefficient}
\label{app:calculation}

We prove \eqref{eq:Apositive} by expanding the original shape energy to
fourth order and eliminating the noncritical variables. The reduction in
Lemma~\ref{lem:reduction} uses the full shape equation.

All quantities in this calculation are Newtonian limits. Set
\[
 q=1,\quad p=3-\sqrt7,\quad r_2=p+1,\quad r_3=0,\quad
 U=\frac{2+\sqrt7}{15},\quad\beta_*=175,\quad L_\ell=\ell(\ell+1).
\]
Along $w=tv$, $v\in V_3$, the radial height is
\begin{equation}\label{eq:height}
 h=tv-t^2v^2+\frac53t^3v^3+O(t^4).
\end{equation}
Add homogeneous corrections of orders $t^j$ to the interior and
exterior radial ball potentials, with traces $A_j,B_j$ at $r=1$.
The radial factors on degree $\ell$ are $r^\ell$ and $r^{-\ell-1}$.
Let $[f]=f^-(1)-f^+(1)$ and let $D_-$ and $D_+$ have multipliers
$\ell$ and $-\ell-1$. The Poisson equation gives
\[
 [u_0'']=-1,\qquad [u_0''']=2-p,\qquad
 [u_0'''']=-8+2p-r_2.
\]
Expanding continuity of the potential and its radial derivative at
$1+h$ yields
\begin{align}\label{eq:transmission}
 A_1=B_1&=v/7,\notag\\
 A_2-B_2&=-v^2/2,&
 D_-A_2-D_+B_2&=pv^2/2,\notag\\
 A_3-B_3&=(2-p)v^3/3,&
 D_-A_3-D_+B_3&=
 \left(-6-\frac p3+\frac{r_2}{6}\right)v^3
 +\frac v2\sum_{\ell=0,2,4,6}L_\ell\Pi_\ell(v^2).
\end{align}
At second order the base solution contributes $-v^2/2$ to the value
jump, and the first-order correction contributes $v^2$. This gives
$A_2-B_2=-v^2/2$. In the third-order
derivative equation one uses
$[u_1''']=18v$ and
$[u_2'']=-pv^2-\frac12\sum L_\ell\Pi_\ell(v^2)$, which follow
from the radial homogeneous equation. These identities fix the signs in \eqref{eq:transmission}.

For the system $A-B=J$, $\ell A+(\ell+1)B=K$,
the exterior trace is $B=(K-\ell J)/(2\ell+1)$. This gives
\begin{align}\label{eq:traces}
 B_2&=\sum_{\ell=0,2,4,6}\frac{p+\ell}{2(2\ell+1)}\Pi_\ell(v^2),\notag\\
 \ip v{B_3}&=\frac{-8+2p/3+r_2/6}{7}\J(v)
 +\frac1{14}\sum_{\ell=0,2,4,6}L_\ell Q_\ell(v).
\end{align}
If $\phi_j$ is the coefficient of the exterior potential evaluated at
the moving boundary, then
\begin{align*}
 \phi_1&=(-U+1/7)v,\\
 \phi_2&=(2U-4/7)v^2+B_2,\\
 \phi_3&=(-14U/3+2)v^3
 -v\sum_{\ell=0,2,4,6}(\ell+1)\Pi_\ell B_2+B_3.
\end{align*}
The density at that boundary is
\[
 \rho(1+h)=1+ptv+(-p+r_2/2)t^2v^2
 +(5p/3-r_2+r_3/6)t^3v^3+O(t^4).
\]
Multiplication gives the nonlocal Euler gradient in exact volume
coordinates, by \eqref{eq:area-gradient}.

\paragraph{Quadratic forcing.}
The cubic perimeter energy is
$P_3(w)=\int(4w^3/3-2w|\nabla w|^2)$, whose gradient at $v\in V_3$
is $\Delta(v^2)-20v^2$. So the quadratic coefficient of the full
projected gradient is
\begin{equation}\label{eq:forcing}
 g_2(v)=\sum_{\ell=2,4,6}g_\ell\Pi_\ell(v^2),\qquad
 g_\ell=-(L_\ell+20)+175\left[
 U\left(2-2p+\frac{r_2}{2}\right)+\frac{p-4}{7}
 +\frac{p+\ell}{2(2\ell+1)}\right].
\end{equation}
The constant channel is removed by the exact volume constraint.

\paragraph{Direct quartic coefficient.}
Pairing the cubic Euler-gradient coefficient with $v$ and dividing by
four gives the direct quartic interaction energy. The part multiplying
$Q_\ell$, modulo a multiple of $\J$, is
\[
 \frac{175}{4}\left[
 \frac{p^2-p-L_\ell}{2(2\ell+1)}+\frac{L_\ell}{14}\right]Q_\ell.
\]
For the perimeter, expansion of the exact area integrand gives
\[
 P_4(v)=\int\left(-\frac73v^4+7v^2|\nabla v|^2-
 \frac18|\nabla v|^4\right)
 =\frac{77}{3}\J-\frac18\sum_{\ell=0,2,4,6}
 \left(12-\frac{L_\ell}{2}\right)^2Q_\ell.
\]
Here $|\nabla v|^2=\frac12\Delta(v^2)+12v^2$ and
$\int v^2|\nabla v|^2=4\int v^4$. We obtain
\begin{equation}\label{eq:direct-coeff}
 Q_{\rm dir}\equiv\sum_{\ell=0,2,4,6}d_\ell Q_\ell
 \pmod{\R\J},\quad
 d_\ell=-\frac18\left(12-\frac{L_\ell}{2}\right)^2
 +\frac{175}{4}\left[
 \frac{p^2-p-L_\ell}{2(2\ell+1)}+\frac{L_\ell}{14}\right].
\end{equation}

\paragraph{Exact complement elimination.}
The quadratic correction is $z_2=-L_W^{-1}g_2$. In degree four it
contributes
$\ip{g_2}{z_2}+\ip{L_Wz_2}{z_2}/2=-\ip{g_2}{L_W^{-1}g_2}/2$.
This gives
\begin{equation}\label{eq:schur-quartic}
 Q=Q_{\rm dir}-\frac12\sum_{\ell=2,4,6}
 \frac{g_\ell^2}{\lambda_\ell}Q_\ell.
\end{equation}
There is a constant channel in \eqref{eq:direct-coeff}, but none in
this elimination sum.

Define $\mathcal L(Q)=10Q(e_T)-9Q(e_S)$. Since $j_S=10j_T/9$,
$\mathcal L(\N)=1$ and $\mathcal L(\J)=0$, so $\mathcal L(Q)$ is
exactly the coefficient $A$ of $\N$. The required data are
\begin{equation}\label{eq:coefficient-table}
\begin{array}{c|c|c|c}
 \ell&\lambda_\ell&g_\ell&\pi\mathcal L(Q_\ell)\\\hline
 0&\text{excluded}&\text{excluded}&1/4\\
 2&4&(674-185\sqrt7)/6&-5/4\\
 4&22/9&(1420-415\sqrt7)/18&69/44\\
 6&240/13&(3854-1565\sqrt7)/78&-25/44
\end{array}
\end{equation}
The last column follows from these projection values:
\[
\begin{array}{c|cccc}
 \ell&0&2&4&6\\\hline
 \pi Q_\ell(e_T)&1/4&0&21/121&200/1573\\
 \pi Q_\ell(e_S)&1/4&5/36&9/484&11575/56628
\end{array}.
\]
To compute them from \eqref{eq:moments}, write the square of the
homogeneous harmonic cubic $P$ as
$P^2=H_6+r^2H_4+r^4H_2+r^6H_0$. Then
$H_0=2\tr A/35$,
$H_2=2(A-\tr A\,\mathrm{Id}/3)_{ij}x_ix_j/7$,
$H_4=(\Delta(P^2)-36r^2H_2-42r^4H_0)/22$, and the remainder is
$H_6$. Restriction to the sphere gives each projection.

Finally, \eqref{eq:direct-coeff}--\eqref{eq:coefficient-table} give
\begin{equation}\label{eq:positive-final}
 \pi A=\sum_{\ell=0,2,4,6}d_\ell\pi\mathcal L(Q_\ell)
 -\frac12\sum_{\ell=2,4,6}\frac{g_\ell^2}{\lambda_\ell}
 \pi\mathcal L(Q_\ell)
 =\frac{66345100\sqrt7-84268315}{1672704}>0.
\end{equation}
For positive $\kappa$, replace the radial powers by the normalized
modified Bessel solutions and the Poisson equation by the Yukawa
equation. Their boundary derivatives at the finitely many degrees
involved converge to the data above. The denominators
$\lambda_2,\lambda_4,\lambda_6$ stay positive. This proves the finite
coefficient continuity used in Proposition~\ref{prop:quartic}.

\paragraph{Independent algebraic check.}
The ancillary program \texttt{check\_transmission.py} reconstructs the
moving-boundary matching equations, the area expansion, and the
Schur correction from polynomial operations over $\mathbb Q(\sqrt7)$.
It does not use fitted numerical values. The program \texttt{check\_angular.py} checks the harmonic projections
and angular Hessians using exact spherical moments. These programs check
the finite algebraic calculations. The reduction and global minimization
arguments are proved in the text.

\paragraph{Use of generative AI.}
ChatGPT (OpenAI) was used for language refinement and consistency checks.
The authors reviewed and verified all content and take full responsibility
for the final manuscript.

\begingroup
\small\raggedright

\endgroup

\end{document}